\documentclass[11pt]{amsart}
\usepackage{amssymb,mathrsfs}
\usepackage{tikz-cd}
\title{Stable Schottky-Jacobi Forms}

\begin{document}

\author{Jae-Hyun Yang}

\address{Department of Mathematics, Inha University,
Incheon 22212, Korea}
\email{jhyang@inha.ac.kr }


\newtheorem{theorem}{Theorem}[section]
\newtheorem{lemma}{Lemma}[section]
\newtheorem{proposition}{Proposition}[section]
\newtheorem{remark}{Remark}[section]
\newtheorem{definition}{Definition}[section]
\newtheorem{conjecture}{Conjecture}[section]

\renewcommand{\theequation}{\thesection.\arabic{equation}}
\renewcommand{\thetheorem}{\thesection.\arabic{theorem}}
\renewcommand{\thelemma}{\thesection.\arabic{lemma}}
\newcommand{\BR}{\mathbb R}
\newcommand{\BQ}{\mathbb Q}
\newcommand{\BT}{\mathbb T}
\newcommand{\BM}{\mathbb M}
\newcommand{\bn}{\bf n}
\def\charf {\mbox{{\text 1}\kern-.24em {\text l}}}
\newcommand{\BC}{\mathbb C}
\newcommand{\BZ}{\mathbb Z}

\thanks{\noindent{Subject Classification:} Primary 14H40, 14H42, 14K25, 32G20\\
\indent Keywords and phrases: stable Jacobi forms, the Schottky problem, the universal Jacobian locus,\\
\indent stable Schottky-Siegel forms, the universal hyperelliptic locus, stable Schottky-Jacobi forms}


\begin{abstract}
In this article, we prove that there do not exist stable Schottky-Jacobi forms for
the universal Jacobian locus and also prove that there exist non-trivial stable Schottky-Jacobi forms for the universal hyperelliptic locus.
\end{abstract}

\maketitle

\newcommand\tr{\triangleright}
\newcommand\al{\alpha}
\newcommand\be{\beta}
\newcommand\g{\gamma}
\newcommand\gh{\Cal G^J}
\newcommand\G{\Gamma}
\newcommand\de{\delta}
\newcommand\e{\epsilon}
\newcommand\lam{\lambda}
\newcommand\z{\zeta}
\newcommand\vth{\vartheta}
\newcommand\vp{\varphi}
\newcommand\om{\omega}
\newcommand\p{\pi}
\newcommand\la{\lambda}
\newcommand\lb{\lbrace}
\newcommand\lk{\lbrack}
\newcommand\rb{\rbrace}
\newcommand\rk{\rbrack}
\newcommand\s{\sigma}
\newcommand\w{\wedge}
\newcommand\fgj{{\frak g}^J}
\newcommand\lrt{\longrightarrow}
\newcommand\lmt{\longmapsto}
\newcommand\lmk{(\lambda,\mu,\kappa)}
\newcommand\Om{\Omega}
\newcommand\ka{\kappa}
\newcommand\ba{\backslash}
\newcommand\ph{\phi}
\newcommand\M{{\Cal M}}
\newcommand\bA{\bold A}
\newcommand\bH{\bold H}
\newcommand\D{\Delta}

\newcommand\Hom{\text{Hom}}
\newcommand\cP{\Cal P}

\newcommand\cH{\Cal H}

\newcommand\pa{\partial}

\newcommand\pis{\pi i \sigma}
\newcommand\sd{\,\,{\vartriangleright}\kern -1.0ex{<}\,}
\newcommand\wt{\widetilde}
\newcommand\fg{\frak g}
\newcommand\fk{\frak k}
\newcommand\fp{\frak p}
\newcommand\fs{\frak s}
\newcommand\fh{\frak h}
\newcommand\Cal{\mathcal}

\newcommand\fn{{\frak n}}
\newcommand\fa{{\frak a}}
\newcommand\fm{{\frak m}}
\newcommand\fq{{\frak q}}
\newcommand\CP{{\mathcal P}_n}
\newcommand\Hnm{{\mathbb H}_n \times {\mathbb C}^{(m,n)}}
\newcommand\BD{\mathbb D}
\newcommand\BH{\mathbb H}
\newcommand\CCF{{\mathcal F}_n}
\newcommand\CM{{\mathcal M}}
\newcommand\Gnm{\Gamma_{n,m}}
\newcommand\Cmn{{\mathbb C}^{(m,n)}}
\newcommand\Yd{{{\partial}\over {\partial Y}}}
\newcommand\Vd{{{\partial}\over {\partial V}}}

\newcommand\Ys{Y^{\ast}}
\newcommand\Vs{V^{\ast}}
\newcommand\LO{L_{\Omega}}
\newcommand\fac{{\frak a}_{\mathbb C}^{\ast}}

\vskip 1mm

\begin{section}{{\bf Introduction}}
\setcounter{equation}{0}

For a positive integer $g$, we let
\begin{equation*}
\BH_g =\left\{ \tau\in\BC^{(g,g)}\,|\ \tau =\,{}^t\tau, \ {\rm Im}\,\tau>0\ \right\}
\end{equation*}
be the Siegel upper half plane of degree $g$
and let
$$Sp(2g,\BR)=\{ M\in \BR^{(2g,2g)}\ \vert \ ^t\!MJ_gM= J_g\ \}$$
be the symplectic group of degree $g$, where $F^{(k,l)}$ denotes
the set of all $k\times l$ matrices with entries in a commutative
ring $F$ for two positive integers $k$ and $l$, $^t\!M$ denotes
the transposed matrix of a matrix $M$ and
$$J_g=\begin{pmatrix} 0&I_g\\
                   -I_g&0\end{pmatrix}.$$
Then $Sp(2g,\BR)$ acts on $\BH_g$ transitively by
\begin{equation*}
M\cdot\tau=(A\tau+B)(C\tau+D)^{-1},
\end{equation*}
where $M=\begin{pmatrix} A&B\\
C&D\end{pmatrix}\in Sp(2g,\BR)$ and $\Om\in \BH_n.$ Let
$$\G_g=Sp(2g,\BZ)=\left\{ \begin{pmatrix} A&B\\
C&D\end{pmatrix}\in Sp(2g,\BR) \,\big| \ A,B,C,D\
\textrm{integral}\ \right\}$$ be the Siegel modular group of
degree $g$. This group acts on $\BH_g$ properly discontinuously.

\vskip 0.35cm
Let ${\Cal A}_g:=\G_g\backslash \BH_g$ be the Siegel modular variety of degree $g$, that is, the moduli space of $g$-dimensional principally polarized abelian varieties,
and let ${\Cal M}_g$ be the the moduli space of projective curves of genus $g$. Then according to Torelli's theorem, the Jacobi mapping
\begin{equation*}
T_g:{\Cal M}_g \lrt {\Cal A}_g
\end{equation*}
defined by
\begin{equation*}
C \longmapsto J(C):= {\rm the\ Jacobian\ of}\ C
\end{equation*}
is injective. The Jacobian locus $J_g:=T_g({\Cal M}_g)$ is a $(3g-3)$-dimensional subvariety of ${\Cal A}_g$

\vskip 0.35cm
The Schottky problem is to characterize the Jacobian locus or its closure ${\bar J}_g$ in ${\Cal A}_g.$ At first this problem had been investigated from
the analytical point of view : to find explicit equations of $J_g$ (or ${\bar J}_g$) in ${\Cal A}_g$ defined by Siegel modular forms on $\BH_g$, for example,
polynomials in the theta constant $\theta \left[ \begin{matrix} \epsilon \\ \delta \end{matrix} \right](\tau,0)$ and their derivatives.
The first result in this direction was due to Friedrich Schottky \cite{Sch} who gave the simple and beautiful equation satisfied by the theta constants of Jacobians of dimension 4.
Much later the fact that this equation characterizes the Jacobian locus $J_4$ was proved by J. Igusa \cite{I} (see also \cite{Fr4}, \cite{G-SM1} and \cite{H-H}). Past decades there has been some progress on the characterization of Jacobians by some mathematicians.

\vskip 0.1cm For two
positive integers $g$ and $h$, we consider the Heisenberg group
$$H_{\BR}^{(g,h)}=\big\{\,(\la,\mu;\ka)\,|\ \la,\mu\in \BR^{(h,g)},\ \kappa\in\BR^{(h,h)},\
\ka+\mu\,^t\la\ \text{symmetric}\ \big\}$$ endowed with the
following multiplication law
$$\big(\la,\mu;\ka\big)\circ \big(\la',\mu';\ka'\big)=\big(\la+\la',\mu+\mu';\ka+\ka'+\la\,^t\mu'-
\mu\,^t\la'\big)$$
with $\big(\la,\mu;\ka\big),\big(\la',\mu';\ka'\big)\in H_{\BR}^{(g,h)}.$
We define the {\it Jacobi group} $G^J$ of degree $g$ and index $h$ that is the semidirect product of
$Sp(2g,\BR)$ and $H_{\BR}^{(g,h)}$
$$G^J=Sp(2g,\BR)\ltimes H_{\BR}^{(g,h)}$$
endowed with the following multiplication law
$$
\big(M,(\lambda,\mu;\kappa)\big)\cdot\big(M',(\lambda',\mu';\kappa'\,)\big)
=\, \big(MM',(\tilde{\lambda}+\lambda',\tilde{\mu}+ \mu';
\kappa+\kappa'+\tilde{\lambda}\,^t\!\mu'
-\tilde{\mu}\,^t\!\lambda'\,)\big)$$ with $M,M'\in Sp(2g,\BR),
(\lambda,\mu;\kappa),\,(\lambda',\mu';\kappa') \in
H_{\BR}^{(g,h)}$ and
$(\tilde{\lambda},\tilde{\mu})=(\lambda,\mu)M'$. Then $G^J$ acts
on $\BH_g\times \BC^{(h,g)}$ transitively by
\begin{equation*}
\big(M,(\lambda,\mu;\kappa)\big)\cdot
(\Om,Z)=\Big(M\cdot\Om,(Z+\lambda \Om+\mu)
(C\Omega+D)^{-1}\Big), \end{equation*}
where $M=\begin{pmatrix} A&B\\
C&D\end{pmatrix} \in Sp(2g,\BR),\ (\lambda,\mu; \kappa)\in
H_{\BR}^{(g,h)}$ and $(\Om,Z)\in \BH_g\times \BC^{(h,g)}.$ We note
that the Jacobi group $G^J$ is {\it not} a reductive Lie group and
the homogeneous space ${\mathbb H}_g\times \BC^{(h,g)}$ is not a
symmetric space. From now on, for brevity we write
$\BH_{g,h}=\BH_g\times \BC^{(h,g)}.$ The homogeneous space $\BH_{g,h}$ is called the
{\it Siegel-Jacobi space} of degree $g$ and index $h$.

\vskip 0.1cm
Let $\G_g^J:=\G_g \ltimes H_\BZ^{(g,h)}$ be the Jacobi modular group. Let
\begin{equation*}
\mathcal A_{g,h}:=\G_g^J\,\backslash \BH_{g,h}
\end{equation*}
be the universal abelian variety. Consider the natural projection map
\begin{equation*}
\pi_{g,h}:\mathcal A_{g,h}\lrt \mathcal A_g.
\end{equation*}
Let
\begin{equation*}
J_{g,h}:=\pi_{g,h}^{-1}(J_g)
\end{equation*}
be the universal Jacobian locus and let
\begin{equation*}
Hyp_{g,h}:=\pi_{g,h}^{-1}(Hyp_g)
\end{equation*}
be the universal hyperelliptic locus, where $Hyp_g$ is the hyperelliptic locus in
$\mathcal A_g.$

\vskip 0.25cm
Let $2\mathcal M$ be a positive definite, even unimodular integral symmetric matrix of
degree $h$. According to Theorem 3.6 in \cite{YJH2}, if $g+h>2k+1$ with a nonnegative integer $k$, the Siegel-Jacobi operator
\begin{equation*}
\Psi_{g,\mathcal M}:J_{k,\M}(\G_g)\lrt J_{k,\M}(\G_{g-1})
\end{equation*}
is an isomorphism (see also Theorem 2.2). Using this fact, we define the notion of stable Jacobi forms of weight $k$ and index $\M$. A Jacobi form $F\in J_{k,\M}(\G_g)$
is said to be a {\it Schottky}-{\it Jacobi form} for $J_{g,h}$\,(resp. $Hyp_{g,h}$)
if it vanishes along $J_{g,h}$\,(resp. $Hyp_{g,h}$). In a natural way, we can define the notion of {\it stable\ Schottky}-{\it Jacobi forms} of index $\M$. For precise definitions, we refer to Definition 2.3 and Definition 4.2.

\vskip 0.35cm
The aim of this paper is to prove the non-existence of stable Schottky-Jacobi forms for the universal Jacobian locus and also to prove that there exist non-trivial stable Schottky-Jacobi forms for the universal hyperelliptic locus.

\vskip 0.25cm
This article is organized as follows. In Section 2, we review some properties of the Siegel-Jacobi operator and the notion of stable Jacobi forms introduced by J.-H. Yang
\cite{YJH5}. In Section 3, we review the notion of stable Schottky-Siegel forms and the works that were done recently by G. Codogni and N. I. Shepherd-Barron\,\cite{Cod, Cod-Sh}.
In Section 4, we introduce the notion of stable Schottky-Jacobi forms and prove the following two theorems.

\begin{theorem}
Let $2\mathcal M$ be a positive definite, even unimodular integral symmetric matrix of degree $h$. Then there do not exist stable Schottky-Jacobi forms of index $\mathcal M$
for the universal Jacobian locus.
\end{theorem}

\vskip 2mm
\begin{theorem}
Let $2\mathcal M$ be a positive definite, even unimodular integral symmetric matrix of degree $h$. Then there exist non-trivial stable Schottky-Jacobi forms of index $\mathcal M$ for the universal hyperelliptic locus.
\end{theorem}

\vskip 2mm
In the final section, we make some comments and present several questions.

\vskip 0.51cm \noindent {\bf Notations:} \ \ We denote by
$\BQ,\,\BR$ and $\BC$ the field of rational numbers, the field of
real numbers and the field of complex numbers respectively. We
denote by $\BZ$ and $\BZ^+$ the ring of integers and the set of
all positive integers respectively. $\BR^+$ denotes the set of all positive real numbers.
$\BZ_+$ and $\BR_+$ denote the set of all nonnegative integers and the set of all nonnegative real numbers respectively. The symbol ``:='' means that the expression on the right is the definition of that on the left. For two positive integers $k$ and $l$, $F^{(k,l)}$ denotes the set
of all $k\times l$ matrices with entries in a commutative ring
$F$. For a square matrix $A\in F^{(k,k)}$ of degree $k$,
$\sigma(A)$ denotes the trace of $A$. For any $M\in F^{(k,l)},\
^t\!M$ denotes the transpose of a matrix $M$. $I_n$ denotes the
identity matrix of degree $n$. We put $i=\sqrt{-1}$.

\vskip 0.1cm

\end{section}

\vskip 5mm


\begin{section}{{\bf Stable Jacobi Forms}}
\setcounter{equation}{0}

\newcommand\POB{ {{\partial}\over {\partial{\overline \Omega}}} }
\newcommand\PZB{ {{\partial}\over {\partial{\overline Z}}} }
\newcommand\PX{ {{\partial}\over{\partial X}} }
\newcommand\PY{ {{\partial}\over {\partial Y}} }
\newcommand\PU{ {{\partial}\over{\partial U}} }
\newcommand\PV{ {{\partial}\over{\partial V}} }
\newcommand\PO{ {{\partial}\over{\partial \Omega}} }
\newcommand\PZ{ {{\partial}\over{\partial Z}} }

\vskip 0.21cm
For a non-negative integer $k$, we denote by $[\G_g,k]$ the vector space of all Siegel modular forms of weight $k$. The Siegel $\Phi$-operator
\begin{equation*}
\Phi_g: [\G_g,k]\lrt [\G_{g-1},k]
\end{equation*}
is an important linear map defined by
\begin{equation}
(\Phi_g f)(\tau): =\lim_{t\lrt\infty}
f \begin{pmatrix} \tau & 0 \\ 0 & it \end{pmatrix}, \qquad f\in [\G_g,k],\ \tau\in \BH_{g-1}.
\end{equation}
H. Maass \cite{M1} proved that if $k$ is even and $k>2g$, then $\Phi_g$ is surjective. E. Freitag \cite{Fr1} proved that if $g>2k$,
then $\Phi_g$ is injective. Using the theory of singular modular forms developed by Freitag \cite{Fr2,Fr5}, he showed the following:
\vskip 2mm
(SO1) $\quad [\G_g,k]=0\quad$ for $g>2k,\ k\not\equiv 0
\,({\rm mod}\,4).$
\vskip 2mm
(SO2) $\quad\Phi_g$ is an isomorphism if $g>2k+1.$

\begin{definition}
A collection $(f_g)_{g\geq 0}$ is called a {\sf stable modular form} of weight $k$ if it satisfies the following conditions (SM1) and (SM2):
\vskip 2mm
{\rm (SM1)} \ \ $f_g\in [\G_g,k]$ for all $g\geq 0$.
\vskip 2mm
{\rm (SM2)} \ \ $\Phi_gf_g=f_{g-1}$ for all $g>0.$
\end{definition}

\vskip 3mm
Let $\rho$ be a rational representation of
$GL(g,\mathbb{C})$ on a finite
dimensional complex vector space
$V_{\rho}.$ Let ${\mathcal M}\in \mathbb R^{(h,h)}$ be a
symmetric
half-integral semi-positive definite matrix of degree $h$.
The canonical automorphic factor
$$ J_{\rho,\mathcal M}: G^J\times \BH_{g,h}\lrt GL(V_\rho)$$
for $G^J$ on $\BH_{g,h}$ is given as follows\,:
\begin{eqnarray*}
 J_{\rho,\mathcal M}((M,(\lambda,\mu;\kappa)),(\tau,z))
&=&e^{2\,\pi\, i\,\sigma\left( {\mathcal M}(z+\lambda\,
\tau+\mu)(C\tau+D)^{-1}C\,{}^t(z+\lambda\,\tau\,+\,\mu)\right) }\\
& &\times\, e^{-2\,\pi\, i\,\sigma\left( {\mathcal M}(\lambda\,
\tau\,{}^t\!\lambda\,+\,2\,\lambda\,{}^t\!z+\,\kappa+
\mu\,{}^t\!\lambda) \right)} \rho(C\,\tau+D),
\end{eqnarray*}
where $M=\left(\begin{matrix} A&B\\ C&D\end{matrix}\right)\in
Sp(2g,\mathbb R),\ (\lambda,\mu;\kappa)\in H_{\mathbb R}^{(g,h)}$
and $(\tau,z)\in \BH_{g,h}.$ We refer to \cite{YJH4} for a geometrical construction of $J_{\rho,\mathcal M}.$

\vskip 0.21cm
Let
$C^{\infty}(\BH_{g,h},V_{\rho})$ be the algebra of all
$C^{\infty}$ functions on $\BH_{g,h}$
with values in $V_{\rho}.$
For $f\in C^{\infty}(\BH_{g,h}, V_{\rho}),$ we define
\begin{eqnarray*}
 \left(f|_{\rho,{\mathcal M}}[(M,(\lambda,\mu;\kappa))]\right)(\tau,z)
&= & J_{\rho,\mathcal M}((M,(\lambda,\mu;\kappa)),(\tau,z))^{-1}\\
& &\ f\left( (A\tau+B)(C\tau+D)^{-1},(z+\lambda\,
\tau+\mu)(C\,\tau+D)^{-1}\right),\nonumber
\end{eqnarray*}
where $M=\left(\begin{matrix} A&B\\ C&D\end{matrix}\right)\in
Sp(2g,\mathbb R),\ (\lambda,\mu;\kappa)\in H_{\mathbb R}^{(g,h)}$
and $(\tau,z)\in \BH_{g,h}.$

\begin{definition}
Let $\rho$ and $\mathcal M$
be as above. Let
$$H_{\mathbb Z}^{(g,h)}:= \left\{ (\lambda,\mu;\kappa)\in
H_{\mathbb R}^{(g,h)}\, \vert
\,\ \lambda,\mu, \kappa\ {\rm integral}\ \right\}$$
be the discrete subgroup of $H_{\mathbb R}^{(g,h)}$.
A $\textsf{Jacobi\ form}$ of index $\mathcal M$
with respect to $\rho$ on a subgroup $\Gamma$ of $\Gamma_g$ of finite index is a holomorphic function $f\in
C^{\infty}(\BH_{g,h},V_{\rho})$ satisfying the following
conditions {\rm (A)} and {\rm (B)}:

\smallskip

\noindent {\rm (A)} \,\ $f|_{\rho,{\mathcal M}}[\tilde{\gamma}] = f$ for
all $\tilde{\gamma}\in {\widetilde\Gamma}:= \Gamma \ltimes
H_{\mathbb Z}^{(g,h)}$.

\smallskip

\noindent {\rm (B)} \,\ For each $M\in\Gamma_g$, $f|_{\rho,\CM}[M]$ has a
Fourier expansion of \\
\indent \ \ \ \ the following form :
$$(f|_{\rho,\CM}[M])(\tau,z) = \sum\limits_{T=\,{}^tT\ge0\atop \text {half-integral}}
\sum\limits_{R\in \mathbb Z^{(g,h)}} c(T,R)\cdot e^{{ {2\pi
i}\over {\lambda_\G}}\,\sigma(T\tau)}\cdot e^{2\pi i\,\sigma(Rz)}$$
with $\lambda_\G(\neq 0)\in\BZ$ and
$c(T,R)\ne 0$ only if $\left(\begin{matrix} { 1\over {\lambda_\G}}T & \frac 12R\\
\frac 12\,^t\!R&{\mathcal M}\end{matrix}\right) \geq 0$.
\end{definition}

\medskip

\indent If $g\geq 2,$ the condition (B) is superfluous by
K{\"o}cher principle\,(\,cf.\,\cite{Zi} Lemma 1.6). We denote by
$J_{\rho,\mathcal M}(\Gamma)$ the vector space of all Jacobi forms
of index $\mathcal{M}$ with respect to $\rho$ on $\Gamma$.
Ziegler\,(\,cf.\,\cite{Zi} Theorem 1.8 or \cite{EZ} Theorem 1.1\,)
proved that the vector space $J_{\rho,\mathcal {M}}(\Gamma)$ is
finite dimensional. In the special case $\rho(A)=(\det(A))^k$ with
$A\in GL(g,\BC)$ and a fixed $k\in\BZ$, we write $J_{k,\CM}(\G)$
instead of $J_{\rho,\CM}(\G)$ and call $k$ the {\it weight} of the
corresponding Jacobi forms. For more results about Jacobi forms with
$g>1$ and $h>1$, we refer to \cite{YJH1,YJH2,YJH3,YJH4,YJH5,YJH6,YJH7} and \cite{Zi}. Jacobi forms play an important role in lifting elliptic cusp forms to Siegel cusp forms of degree $2g$ (cf.\,\cite{Ik, Ik1}).

\vskip 2mm Now we consider the special case $\rho=\det^k$ with $k\in \BZ_+$. We define the {\sf Siegel}-{\sf Jacobi\ operator}
\begin{equation*}
\Psi_{g,\M}: J_{k,\CM}(\G_g)\lrt J_{k,\CM}(\G_{g-1})
\end{equation*}
by
\begin{equation}
(\Psi_{g,\M}F)(\tau,z):= \lim_{t\lrt\infty}
F\left( \begin{pmatrix} \tau & 0 \\ 0 & it \end{pmatrix},(z,0)\right),\\
\end{equation}
where $F \in J_{k,\CM}(\G_g),\ \tau\in \BH_{g-1}$ and
$z\in \BC^{(h,g-1)}.$ We observe that the above limit exists
and $\Psi_{g,\M}$ is a well-defined linear map\,(cf.\,\cite{Zi}).
\vskip 0.1cm
J.-H. Yang\,\cite{YJH2} proved the following theorems.
\begin{theorem}
Let $2\M$ be a positive even unimodular symmetric integral matrix of degree $h$ and let $k$ be an even nonnegative integer. If $g+h>2k$, then the Siegel-Jacobi operator
$\Psi_{g,\M}$ is injective.
\end{theorem}
\vskip 1mm\noindent
{\it Proof.} See Theorem 3.5 in \cite{YJH2}. \hfill $\Box$

\begin{theorem}
Let $2\M$ be as above in Theorem 2.1 and let $k$ be an even nonnegative integer. If $g+h>2k+1$, then the Siegel-Jacobi operator $\Psi_{g,\M}$ is an isomorphism.
\end{theorem}
\vskip 1mm\noindent
{\it Proof.} See Theorem 3.6 in \cite{YJH2}. \hfill $\Box$

\begin{remark}
A Jacobi form in $J_{k,\M}(\G_g)$ is said to be {\sf singular} if it admits a Fourier expansion such that a Fourier coefficient $c(T,R)$
vanishes unless
\begin{equation*}
\det \begin{pmatrix}
       T & {\frac 12}R \\
       {\frac 12}\,{}^t R & \M
     \end{pmatrix}=0.
\end{equation*}
Let $2\M$ be as above in Theorem 2.1. Yang proved that if $k$ is an even nonnegative integer and $g+{\rm rank}(\M)>2k$, then any non-zero Jacobi form in $J_{k,\M}(\G_g)$ is singular
(cf.\,\cite[Theorem 4.5]{YJH3}).
\end{remark}

\begin{theorem}
Let $2\M$ be as above in Theorem 2.1 and let $k$ be an even nonnegative integer. If $2k>4g+h$,
then the Siegel-Jacobi operator $\Psi_{g,\M}$ is surjective.
\end{theorem}
\vskip 1mm\noindent
{\it Proof.} See Theorem 3.7 in \cite{YJH2}. \hfill $\Box$

\begin{remark}
Yang \cite[Theorem 4.2]{YJH2} proved that the action of the Hecke operatos on Jacobi forms is compatible with that of the Siegel-Jacobi operator.
\end{remark}

\begin{definition}
A collection $(F_g)_{g\geq 0}$ is called a {\sf stable Jacobi form} of weight $k$ and index $\M$ if it satisfies the following conditions (SJ1) and (SJ2):
\vskip 2mm
{\rm (SJ1)} \ \ $F_g\in J_{k,\CM}(\G_g)$ \ \ for all $g\geq 0.$
\vskip 2mm
{\rm (SJ2)} \ \ $\Psi_{g,\M}F_g=F_{g-1}$ \ \ for all $g\geq 1.$
\end{definition}

\begin{remark}
The concept of a stable Jacobi forms was introduced by Yang\,\cite{YJH5}.
\end{remark}

\noindent {\bf Example.}
Let $S$ be a positive even unimodular symmetric integral matrix of degree $2k$ and let $c\in \BZ^{(2k,h)}$ be an integral matrix. We define the theta series $\vartheta_{S,c}^{(g)}$ by
\begin{equation*}
\vartheta_{S,c}^{(g)}(\tau,z):=\sum_{\la\in\BZ^{(2k,g)}}
e^{\pi i\left\{ \sigma(S\la\tau\,{}^t\la)+ 2\sigma({}^tc\,S\la\,{}^tz)\right\} }, \quad (\tau,z)\in \BH_{g,h}.
\end{equation*}
It is easily seen that $\vartheta_{S,c}^{(g)}\in J_{k,\CM}(\G_g)$
with $\M:={\frac 12}~\!{}^t\!cSc$ for all $g\geq 0$ and $\Psi_{g,\M}\vartheta_{S,c}^{(g)}=\vartheta_{S,c}^{(g-1)}$ for all
$g\geq 1.$ Thus the collection
\begin{equation*}
\Theta_{S,c}:=\left( \vartheta_{S,c}^{(g)} \right)_{g\geq 0}
\end{equation*}
is a stable Jacobi form of weight $k$ and index $\M$.

\end{section}

\vskip 0.51cm


\begin{section}{{\bf Stable Schottky-Siegel Forms}}
\setcounter{equation}{0}

\vskip 0.21cm

Let $\mathcal A_g^{\rm Sat}$ be the Satake compactification of the Siegel modular variety $\mathcal A_g$ \,(cf.\,\cite{Sat}).
\begin{equation*}
\mathcal A_g^{\rm Sat}= \mathcal A_g \cup \mathcal A_{g-1}\cup\cdots
\mathcal A_1\cup \mathcal A_0.
\end{equation*}

W. Baily \cite{B} proved that $\mathcal A_g^{\rm Sat}$ is a normal
projective variety in which $\mathcal A_g$ is Zariski open. In particular, we have a closed embedding
\begin{equation*}
\iota_g:\mathcal A_{g-1}^{\rm Sat} \hookrightarrow \mathcal A_{g}^{\rm Sat}.
\end{equation*}
The collection $(\mathcal A_{g}^{\rm Sat})_{g\geq 0}$ and the above embeddings $(\iota_g)_{g\geq 0}$ define the projective limit
\begin{equation*}
\mathcal A_\infty^{\rm Sat}:=\bigcup_{g\geq 0}
\mathcal A_g^{\rm Sat}=
\lim_{\begin{subarray}{c} \longleftarrow\\ ^g \end{subarray}} \mathcal A_g^{\rm Sat}
\end{equation*}
which is called the {\it stable Satake compactification}. Let
$\mathcal L_g$ be the determinant line bundle of the Hodge bundle over $\mathcal A_g$. The we have the isomorphism
\begin{equation*}
H^0(\mathcal A_g,\mathcal L_g^{\otimes k}) \cong [\G_g,k].
\end{equation*}
Let $J_g^{\rm Sat}$ (resp. $Hyp_g^{\rm Sat}$) be the closure of $J_g$ (resp. $Hyp_g$) inside $\mathcal A_g^{\rm Sat}$. We define
\begin{equation*}
J_\infty:=\bigcup_{g\geq 0} J_g^{\rm Sat}\quad {\rm and}\quad
Hyp_\infty:=\bigcup_{g\geq 0} Hyp_g^{\rm Sat}.
\end{equation*}

\begin{definition}
A pair $(\Lambda,Q)$ is called a quadratic form if $\Lambda$ is a lattice and $Q$ is an integer-valued bilinear symmetric form on
$\Lambda$. The rank of $(\Lambda,Q)$ is defined to be the rank of
$\Lambda$. For $v\in \Lambda$, the integer $Q(v,v)$ is called the norm of $v$. A quadratic form $(\Lambda,Q)$ is said to be even if
$Q(v,v)$ is even for all $v\in \Lambda$.
A quadratic form $(\Lambda,Q)$ is said to be unimodular if $\det (Q)=1.$
\end{definition}

\begin{definition}
Let $(\Lambda, Q)$ be an even unimodular positive definite quadratic form of rank $m$.
For a positive integer $g$, the theta series $\theta_{Q,g}$ associated to $(\Lambda, Q)$ is defined to be
\begin{equation*}
\theta_{Q,g}(\tau):=\sum_{x_1,\cdots,x_g\in \Lambda}\exp \left( \pi i \sum_{p,q=1}^g Q(x_p,x_q)\tau_{pq}\right),\qquad \tau=(\tau_{pq})\in \BH_g.
\end{equation*}
\end{definition}

It is well known that $\theta_{Q,g}(\tau)$ is a Siegel modular form on $\BH_g$ of weight ${\frac m2}.$ We easily see that
\begin{equation*}
\Phi_{g+1}(\theta_{Q,g+1})=\theta_{Q,g}\qquad {\rm for\ all}\ g \geq 0.
\end{equation*}
Therefore the collection of all theta series associated to $(\Lambda, Q)$
\begin{equation}
\Theta_Q:=\left( \theta_{Q,g} \right)_{g\geq 0}
\end{equation}
is a stable modular form of weight ${\frac m2}.$

\vskip 3mm
Freitag \cite{Fr2} proved the following theorem.
\begin{theorem}
The ring of stable modular forms is a polynomial ring in countably many
theta series $\Theta_{Q}=\left( \theta_{Q,g} \right)_{g\geq 0}$ associated to irreducible positive even
unimodular quadratic forms.
\end{theorem}
\vskip 1mm\noindent
{\it Proof.} See Theorem 2.5 in \cite{Fr2}. \hfill $\Box$

\begin{definition}
A modular form $f\in [\G_g,k]$ is called a {\sf Schottky-Siegel form} of weight $k$ for $J_g$ (resp. $Hyp_g$) if it vanishes along
$J_g$ (resp. $Hyp_g$). A collection $(f_g)_{g\geq 0}$ is called a {\sf stable Schottky-Siegel form} of weight $k$ for the Jacobian locus (resp. the hyperelliptic locus) if $(f_g)_{g\geq 0}$ is a stable modular form of weight $k$ and $f_g$ vanishes along $J_g$ (resp. $Hyp_g$) for every $g\geq 0.$
\end{definition}

G.~Codogni and N.~I. Shepherd-Barron \cite{Cod-Sh} proved the following.
\begin{theorem}
There do not exist stable Schottky-Siegel form for the Jacobian locus.
\end{theorem}
\vskip 1mm\noindent
{\it Proof.} See Theorem 1.3 and Corollary 1.4 in \cite{Cod-Sh}. \hfill $\Box$

\begin{remark}
Let
\begin{equation}
  \varphi_g(\tau):=\theta_{E_8\oplus E_8,g}(\tau)
  -\theta_{D_{16}^+,g}(\tau),\qquad \tau\in\BH_g
\end{equation}
be the Igusa modular form, that is, the difference of the theta series in genus $g$ associated to the two distinct positive even unimodular quadratic forms $E_8\oplus E_8$ and $D_{16}^+$ of rank $16$. We see that $\varphi_g(\tau)$ is a Siegel modular form on $\BH_g$ of weight $8$. Since $\Phi_g \varphi_g=\varphi_{g-1}$ for all $g\geq 1$, a collection $(\varphi_g)_{g\geq 0}$ is a stable
modular form of weight $8$. Igusa \cite{I,I1} showed that the Schottky-Siegel form discovered by Schottky \cite{Sch} is an explicit rational multiple of $\varphi_4$. In \cite{I}, he also
showed that the Jacobian locus $J_4$ is reduced and irreducible, and so cuts out exactly $J_4$ in $\mathcal A_4.$ Indeed, $\varphi_4(\tau)$ is a degree $16$ polynomial in the Thetanullwerte of genus $4$. On the other hand, Grushevsky and Salvati Manni \cite{G-SM2} showed that the Igusa modular form $\varphi_5$ of genus $5$ cuts out exactly the trigonal locus in $J_5$ and so does not vanish along $J_5$. Thus $(\varphi_g)_{g\geq 0}$ is not a stable
Schottky-Siegel form.
\end{remark}

\vskip 2mm
G.~Codogni \cite{Cod} proved the following.
\begin{theorem}
There exist non-trivial stable Schottky-Siegel form for the hyperelliptic locus. Precisely the ideal of stable Schottky-Siegel forms for the hyperelliptic locus is generated by differences of theta series
\begin{equation*}
\Theta_P-\Theta_Q,
\end{equation*}
where $P$ and $Q$ are positive definite even unimodular quadratic forms of the same rank.
\end{theorem}
\vskip 1mm\noindent
{\it Proof.} See Theorem 1.2 in \cite{Cod}. \hfill $\Box$

\begin{remark}
Let $P$ and $Q$ be two positive even
unimodular quadratic forms of the same rank. We let
\begin{equation*}
\Theta_P:=\left( \theta_{P,g} \right)_{g\geq 0}\quad
{\rm and}\quad
\Theta_Q:=\left( \theta_{Q,g} \right)_{g\geq 0}
\end{equation*}
be two stable modular forms. Codogni \cite[Theorem 1.4]{Cod} showed that
the difference of theta series
\begin{equation*}
\Theta_P-\Theta_Q
\end{equation*}
is a stable Schottky-Siegel form for the hyperelliptic locus when
one of the following conditions {\rm (1)--(3)}:
\vskip 2mm
{\rm (1)} {\rm rank}($P$)={\rm rank}($Q$)=24\ and the two quadratic forms have the same number of vectors\\
\indent\ \ \ \  of norm $2$;
\vskip 2mm
{\rm (2)} {\rm rank}($P$)={\rm rank}($Q$)=32\ and the two quadratic forms do not have any vectors of norm \\
\indent\ \ \ \ \ $2$;
\vskip 2mm
{\rm (3)} {\rm rank}($P$)={\rm rank}($Q$)=48\ and the two quadratic forms do not have any vectors of norm \\
\indent\ \ \ \ \ $2$ and $4$.
\end{remark}

\end{section}

\vskip 0.21cm


\begin{section}{{\bf Stable Schottky-Jacobi Forms and Proofs of Main Theorems}}
\setcounter{equation}{0}
\vskip 0.21cm
In this section, we introduce the notion of stable Schottky-Jacobi forms and prove the main theorems.
\vskip 1mm
We let
\begin{equation*}
\mathcal A_{g,h}:=\G^J_{g,h}\,\backslash \BH_{g,h}
\end{equation*}
be the universal abelian variety and let
\begin{equation*}
\mathcal A_{g,h}^{\rm Sat}:=\mathcal A_{g,h}\cup \mathcal A_{g-1,h}
\cup\cdots \cup \mathcal A_{1,h} \cup \mathcal A_{0,h}
\end{equation*}
be the Satake compactification of $\mathcal A_{g,h}$. We consider the natural projection map
\begin{equation*}
\pi_{g,h}:\mathcal A_{g,h}\lrt \mathcal A_{g}
\end{equation*}
of $\mathcal A_{g,h}$ onto $\mathcal A_{g}.$ Let
\begin{equation*}
J_{g,h}:=\pi_{g,h}^{-1}(J_g)
\end{equation*}
be the universal Jacobian locus and let
\begin{equation*}
Hyp_{g,h}:=\pi_{g,h}^{-1}(Hyp_g)
\end{equation*}
be the universal hyperelliptic locus. Let ${J}_{g,h}^S$
(resp. $Hyp_{g,h}^S$) be the closure of $J_{g,h}$ (resp. $Hyp_{g,h}$) in $\mathcal A_{g,h}^{\rm Sat}$.
We put
\begin{equation*}
\mathcal A_{\infty,h}:=\bigcup_{g\geq 0}
\mathcal A_{g,h}^{\rm Sat},
\end{equation*}

\begin{equation*}
 J_{\infty,h}:=\bigcup_{g\geq 0} J_{g,h}^{\rm S}
\end{equation*}
and
\begin{equation*}
 Hyp_{\infty,h}:=\bigcup_{g\geq 0} Hyp_{g,h}^{\rm S}.
\end{equation*}

\begin{definition}
Let $\M$ be a half-integral semi-positive symmetric matrix of degree $h$ and $k\in\BZ_+.$ A Jacobi form $F\in J_{k,\CM}(\G_g)$ is called a {\sf Schottky-Jacobi form} of weight $k$ and index $\M$ for the universal Jacobian (resp. hyperelliptic) locus if it vanishes along $J_{g,h}$ (resp. $Hyp_{g,h}$).
\end{definition}

\begin{definition}
Let $\M$ be a half-integral semi-positive symmetric matrix of degree $h$ and $k\in\BZ_+.$
A collection $(F_g)_{g\geq 0}$ is called a {\sf stable Schottky-Jacobi form} of weight $k$ and index $\M$ if it satisfies the following conditions {\rm (SSJ1)} and {\rm (SSJ2)}:
\vskip 2mm
{\rm (SSJ1)} \ \ $F_g\in J_{k,\M}(\G_g)$ is a Schottky-Jacobi form of weight $k$ and index $\M$ for all $g\geq 0.$
\vskip 2mm
{\rm (SSJ2)}\ \ $\Psi_{g,\M}F_g=F_{g-1}$ for all $g\geq 1.$
\end{definition}

\begin{theorem}
Let $2\M$ be a positive even unimodular symmetric integral matrix of degree $h$. Then there do not exist stable Schottky-Jacobi forms of index $\M$ for the universal Jacobian locus.
\end{theorem}
\vskip 1mm\noindent
{\it Proof.} We first observe that $h\equiv 0\,({\rm mod}\,8)$ (cf.\,\cite{Ser}). Assume that there exists a non-trivial stable Schottky-Jacobi form $(F_g)_{g\geq 0}$ of weight $k$ and index $\M$ for the universal Jacobian locus.

\vskip 2mm\noindent
{\bf Case 1:} $k$ is even.
\vskip 1mm
Using the Shimura isomorphsm\,(cf.\,\cite{Sh1, Sh2}), we obtain the following
\begin{equation}
J_{k,\M}(\G_g)=[\G_g,k_*]\cdot \vartheta_{2\M}^{[g]}(\tau,z),
\end{equation}
where $k_*:=k-{h\over 2}$ and
\begin{equation}
\vartheta_{2\M}^{[g]}(\tau,z):=\sum_{\la\in\BZ^{(h,g)}}
e^{2\pi i\,\sigma\left( \M (\la \tau\,{}^t\la+2\la\,{}^tz)\right)}.
\end{equation}
We refer to \cite[Theorem 3.3]{Zi} for the proof of the formula (4.1). We see from (SO1) in Section 2 that $[\G_g,k_*]=0$ if $g+h>2k$ and
$k\not\equiv 0\,({\rm mod}\ 4)$. So $k\equiv 0 \,({\rm mod}\ 4)$.
We observe that the Siegel-Jacobi operator
$\Psi_{g,\M}:J_{k,\M}(\G_g)\lrt J_{k,\M}(\G_{g-1})$ is an isomorphism if $g+h>2k+1$ (see Theorem 2.2 in Section 2).
It is easy to see that
\begin{equation*}
\Psi_{g,\M}\vartheta_{2\M}^{[g]}=\vartheta_{2\M}^{[g-1]}
\quad {\rm for\ all}\ g\geq 1.
\end{equation*}
According to the formula (4.1), we may write
\begin{equation*}
F_g(\tau,z)=f_g(\tau)\cdot \vartheta_{2\M}^{[g]}(\tau,z),\quad f\in [\G,k_*].
\end{equation*}
Now we have, for $(\tau,z)\in\BH_{g-1,h},$
\begin{eqnarray*}
\left( \Psi_{g,\M}F_g\right)(\tau,z) &=&
\lim_{t\lrt\infty}F_g\left(
\begin{pmatrix}
  \tau & 0 \\
  0 & it
\end{pmatrix},(z,0)\right)\\
&=&
\lim_{t\lrt\infty}f_g
\begin{pmatrix}
  \tau & 0 \\
  0 & it
\end{pmatrix}\cdot
\vartheta_{2\M}^{[g]}
\left(
\begin{pmatrix}
  \tau & 0 \\
  0 & it
\end{pmatrix},(z,0)\right)\\
&=&
\left( \Phi_{g}f_g\right)(\tau)\cdot \vartheta_{2\M}^{[g-1]}(\tau,z).
\end{eqnarray*}
Here $\Phi_g$ is the Siegel $\Phi$-operator defined by (2.1).
\vskip 1mm
On the other hand, by the assumption that $(F_g)_{g\geq 0}$ is a stable Schottky-Jacobi form, we have
\begin{equation*}
  \Psi_{g,\M}F_g=F_{g-1}=f_{g-1}\cdot \vartheta_{2\M}^{[g-1]}
  \quad {\rm for\ some}\ f_{g-1}\in [\G_{g-1},k]
\end{equation*}
for all $g\geq 1.$
Therefore
\begin{equation*}
  \Phi_g f_g=f_{g-1}\quad {\rm for\ all}\ g\geq 1.
\end{equation*}
Obviously $f_g$ vanishes along $J_g$ for all $g\geq 0$. Thus
$(f_g)_{g\geq 0}$ is a non-trivial stable Schottky-Siegel form of weight $k_*$. This contradicts the non-existence of
a non-trivial stable Schottky-Siegel form for the Jacobian locus
(see Theorem 3.2).

\vskip 3mm\noindent
{\bf Case 2:} $k$ is odd.
\vskip 1mm
Using the Shimura isomorphism, we may write
\begin{equation*}
  F_g(\tau,z)=\psi_g(\tau) \cdot \vartheta_{2\M}^{[g]}(\tau,z)
  \qquad {\rm for\ all}\ g\geq 1,
\end{equation*}
where $\vartheta_{2\M}^{[g]}(\tau,z)$ is the theta series defined by Formula (4.2) and $f_g(\tau)$ satisfies the following behaviours (4.3) and (4.4):
\begin{equation}
\psi_g(\tau+S)=\psi_g(\tau)\qquad {\rm for\ all}\ S=\,{}^t\!S\in \BZ^{(g,g)};
\end{equation}

\begin{equation}
\psi_g(-\tau^{-1})=\det (-\tau)^k \det \left( {\frac{\tau}{i}} \right)^{-{\frac h2}} \psi_g(\tau),\quad \tau\in\BH_g.
\end{equation}
We put
\begin{equation*}
\xi_g(\tau):=\left\{ \psi_g(\tau) \right\}^{2}\qquad
{\rm for\ all}\ g\geq 1.
\end{equation*}
Then we see easily that a collection $\xi=(\xi_g)_{g\geq 0}$ is a non-trivial stable Schottky-Siegel form of weight $2k-h$ for the Jacobian locus. This contradicts the non-existence of
a non-trivial stable Schottky-Siegel form for the Jacobian locus.
Hence we complete the proof of the above theorem(=Theorem 1.1).
\hfill $\Box$

\begin{theorem}
Let $2\M$ be a positive even unimodular symmetric integral matrix of degree $h$. Then there exist non-trivial stable Schottky-Jacobi forms of $\M$ for the universal hyperelliptic locus.
\end{theorem}
\vskip 1mm\noindent
{\it Proof.}
According to Theorem 3.3, there exists a non-tivial stable Schottky-Siegel form $(f_g)_{g\geq 0}$ of weight $k$ for the hyperelliptic locus. We see from (SO1) in Section 2 that
$k\equiv 0 \,({\rm mod}\ 4)$ and $k\in\BZ^+.$
We put $\ell:=k+{h\over 2}.$ Then using the Shimura isomorphism,
we have
\begin{equation*}
J_{\ell,\M}(\G_g)=[\G_g,k]\cdot \vartheta_{2\M}^{[g]}(\tau,z),
\end{equation*}
where $\vartheta_{2\M}^{[g]}(\tau,z)$ is the theta series defined by Formula (4.2). We define the Jacobi forms
\begin{equation*}
F_g(\tau,z):=f_g(\tau)\cdot \vartheta_{2\M}^{[g]}(\tau,z),
\quad g\geq 0.
\end{equation*}
Since $f_g\in [\G_g,k]$ is a Jacobi form of weight $k$ and index $0$, we get $F_g\in J_{\ell,\M}(\G_g)$ for all $g\geq 0.$ For
$[(\tau,z])\in Hyp_{g,h}$,
\begin{equation*}
F_g(\tau,z)=0,\quad g\geq 0.
\end{equation*}
By a simple calculation, we obtain
\begin{equation*}
\Psi_{g,\M}F_g=F_{g-1}\qquad {\rm for\ all}\ g\geq 1.
\end{equation*}
Thus
$(F_g)_{g\geq 0}$ is a non-trivial stable Schottky-Jacobi form
of weight $\ell$ and index $\M$ for the universal hyperelliptic
locus $Hyp_{\infty,h}$. This completes the proof of the above theorem(=Theorem 1.2).
\hfill $\Box$

\vskip 2mm
We define the invariant $\mu(Q)$ of a quadratic form $(Q,\Lambda)$
by
\begin{equation*}
\mu(Q):={\rm min}\{ Q(v,v)\,|\ v\in \Lambda, v\neq 0\,\}.
\end{equation*}

\begin{theorem}
Let $2\M$ be a positive even unimodular symmetric integral matrix of degree $h$. Let $(Q,\Lambda)$ and $(P,\G)$ be two positive even unimodular quadratic forms of rank $m$. Assume that
\begin{equation*}
{m\over \mu}\leq 8,\qquad where\ \mu:={\rm min}\{ \mu (Q),\mu(P)\}.
\end{equation*}
We put
\begin{equation*}
F_g(\tau,z):=\left\{ \theta_{Q,g}(\tau)-\theta_{P,g}(\tau)\right\}
\cdot \vartheta_{2\M}^{(g)}(\tau,z),\quad g\geq 0.
\end{equation*}
Then $(F_g)_{g\geq 0}$ is a stable Schottky-Jacobi form of weight ${\frac 12}(m+h)$ and index $\M$ for the universal hyperelliptic locus.
\end{theorem}
\vskip 1mm\noindent
{\it Proof.}
It is easily seen that
\begin{equation}
\theta_{Q,g},\,\theta_{P,g}\in J_{{m\over 2},0}(\G_g) \quad {\rm and} \quad \theta_{Q,g}\cdot\vartheta_{2\M}^{[g]},\
\theta_{P,g}\cdot\vartheta_{2\M}^{[g]}\in
J_{{1\over 2}(m+h),\M}(\G_g)
\end{equation}
for all $g\geq 0$. The proof follows immediately from Theorem 5.5 in \cite{Cod} and the above facts (4.5).
\hfill $\Box$

\end{section}

\vskip 5.3mm
\begin{section}{{\bf Final Remarks}}
\setcounter{equation}{0}
\vskip 0.21cm
In the final section, we make some remarks and present several open questions.

\begin{remark}
Let $2\M$ be a positive even unimodular symmetric integral matrix of degree $h$. Assume that
\begin{equation*}
g+{\rm rank}(\M)> 2k+1 \quad {\rm and} \ k\in\BZ_+\ {\rm is\ even}.
\end{equation*}
We denote by ${\mathcal C}_{k,\M}$ be the vector space of stable Jacobi forms of weight $k$ and index $\M$. According to Theorem 2.2,
the Siegel-Jacobi operator
$\Psi_{g,\M}:J_{k,\M}(\G_g)\lrt J_{k,\M}(\G_{g-1})$
is an isomorphism, and hence we obtain
\begin{equation*}
\dim {\mathcal C}_{k,\M}=\dim J_{k,\M}(\G_g).
\end{equation*}
From Formula (4.1), we see that
\begin{equation*}
J_{k,\M}(\G_g)=[\G_g,k_*]\cdot \vartheta_{2\M}^{[g]}(\tau,z),\qquad {\rm where}\ k_*:=k-{\frac h2}.
\end{equation*}
Therefore from (SO1) in Section 2, we get the vanishing result:
\begin{equation*}
J_{k,\M}(\G_g)=0 \qquad {\rm if}\ 2k\not\equiv h\ ({\rm mod}\,8).
\end{equation*}
Thus $k\equiv 0\,({\rm mod}\,4)$ if $J_{k,\M}(\G_g)\neq 0.$
According to Yang \cite{YJH3}, any Jacobi form in $J_{k,\M}(\G_g)$
is {\sf singular}. We note that any element in $[\G_g,k-{\frac h2}]$
is a singular modular form (see \cite{Fr2, Fr5}). Hence we conclude that ${\mathcal C}_{k,\M}$ is spanned by stable Jacobi forms of the form
 \begin{equation*}
\left( \theta_{P,g}(\tau)\vartheta_{2\M}^{[g]}(\tau,z)
\right)_{g \geq 0},
\end{equation*}
where $P$ runs over the set of positive even unimodular quadratic forms of rank $2k-h.$
\end{remark}

\begin{remark}
Let $\varphi_g(\tau)$ be the Igusa modular form defined by the formula (3.2). We denote by $[\G_g,k]_0$ be the space of all Siegel cuspidal Hecke eigenforms on $\BH_g$ of weight $k$.
It is known that $[\G_4,8]_0=\BC\cdot\varphi_4$ (for a nice proof of this, we refer to \cite{DI}).
Poor \cite{P} showed that $\varphi_g(\tau)$ vanishes on the hyperelliptic locus $Hyp_g$ for all $g\geq 1$, and the divisor of $\varphi_g(\tau)$ in $\mathcal{A}_g$ is proper and irreducible for all $g\geq 4$.
And Ikeda \cite{Ik, Ik1}
proved that if $g\equiv k\,({\rm mod}\,2),$ there exists a canonical lifting
\begin{equation*}
  I_{g,k}:[\G_1,2k]_0\lrt [\G_{2g},g+k]_0.
\end{equation*}
Considering the special cases of the Ikeda lift maps $I_{2,6}$ and $I_{6,6}$, Breulman and Kuss \cite{BK} showed that
\begin{equation*}
  I_{2,6}(\Delta)=c\,\varphi_4, \quad c(\neq 0)\in\BC,
\end{equation*}
and constructed a nonzero Siegel cusp form of degree $12$ and weight $12$ which is the image of $\Delta(\tau)$, where
\begin{equation*}
\Delta(\tau)=(2\pi i)^{12}~q\prod_{n=1}^{\infty}
\left( 1-q^n\right)^{24},\quad q:=e^{2\pi i\tau},\quad \tau\in\BH_1
\end{equation*}
is a cusp form of weight 12.
\end{remark}

\begin{remark}
We consider a half-integral semi-positive symmetric integral
matrix $\M$ such that $2\M$ is {\it not} even or which is {\it not} unimodular.
The natural questions arise:
\vskip 2mm\noindent
{\bf Question\ 1.} Are there non-trivial stable Schottky-Jacobi forms of index $\M$ for the universal Jacobian locus ?
\vskip 2mm\noindent
{\bf Question\ 2.} Are there non-trivial stable Schottky-Jacobi forms of index $\M$ for the universal hyperelliptic locus ?
\end{remark}

\begin{remark}
Let $2\M$ be a positive even unimodular symmetric integral matrix of degree $h$. For two nonnegative integers $k$ and $\ell$, we let $A_{k,\ell\M}$ be the vector space of stable Jacobi forms of weight $k$ and index $\ell\M$. We put
\begin{equation*}
  A(\M):=\bigoplus_{\ell =0}^{\infty}\bigoplus_{k=0}^{\infty}
A_{k,\ell\M}.
\end{equation*}
Then we see easily that
\begin{equation*}
  A_{k,\ell\M}{\bullet} A_{p,q\M}\subset A_{k+p,(\ell+q)\M}
\end{equation*}
with respect to the natural multiplication $\bullet$.
Thus $A(\M)$ is a bi-graded ring.
Let $I(\M)$ be the space of all stable Schottky-Jacobi forms
for the universal hyperelliptic locus contained in $A(\M)$.
Then $I(M)$ is an ideal of $A(\M)$.

\vskip 2mm
According to Theorem 3.1, the subset
\begin{equation*}
  A(\M)_0:=\bigoplus_{k=0}^{\infty}A_{k,0}
\end{equation*}
of $A(\M)$ has a polynomial ring structure.

\vskip 2mm\noindent
Let
\begin{equation*}
  A^{[4]}(\M)_1:=\bigoplus_{k \equiv 0\,({\rm mod}\,4)}A_{k,\M}
\end{equation*}
and let $B^{[4]}(\M)_1$ be the subspace of all stable Schottky-Jacobi forms for the universal hyperelliptic locus contained in $A^{[4]}(\M)_1$. Using Theorem 4.2 in \cite{Cod},
we can show that
\begin{equation}
\left( \Theta_P-\Theta_Q\right)\,\Theta_{2\M}
\end{equation}
is a stable Schottky-Jacobi form for the universal hyperelliptic locus of weight ${\frac 12}(m+h)$ and index $\M$, that is,
$\left( \Theta_P-\Theta_Q\right)\,\Theta_{2\M}\in B^{[4]}(\M)_1$.
Here $P$ and $Q$ are two positive even unimodular quadratic forms of the same rank $m\ (m\in\BZ^+)$,  and $\Theta_P,\ \Theta_Q$ are stable modular forms that are defined in Formula (3.1). The subspace
$B^{[4]}(\M)_1$ of $A^{[4]}(\M)_1$ is spanned by all the stable Jacobi forms of type (5.1), where $m$ runs over the set of all positive integers $8n\,(n\in\BZ^+)$.

\vskip 2mm\noindent
{\bf Question\ 3.} What kinds of structures does $A(\M)$ have ?
\end{remark}

\end{section}


\end{document}